%


\input amstex
\documentstyle{amsppt}
\magnification=\magstep1
\NoBlackBoxes
\hsize 6.0truein
\vsize 8.2truein
\hcorrection{.2truein}
\vcorrection{.2truein}
\loadbold

\define\sg{\Sigma^*_n}
\redefine\sgn{\Sigma^*_{n-1}}
\define\so{\Sigma^*_0}

\define\cof{\operatorname{cof}}
\define\Seq{\operatorname{Seq}}
\define\oM{\overline{M}}
\define\ORD{\operatorname{ORD}}

\redefine\phi{\varphi}

\define\Lim{\operatorname{Lim}}

\define\Dom{\operatorname{Dom}}

\document
\baselineskip=15pt

\font\bigtenrm=cmr12 scaled\magstep2
\centerline
{\bigtenrm {Jensen's $\sum^*$ Theory and the Combinatorial}}

\centerline
{\bigtenrm {Content of $V=L$ }}
\vskip20pt

\font\bigtenrm=cmr10 scaled\magstep2
\centerline
{\bigtenrm{Sy D. Friedman}\footnote"*"{Research supported by NSF contract \#
8903380.}} 

\centerline
{\bigtenrm {M.I.T.}}

\vskip20pt

\comment
lasteqno 1@0
\endcomment

An awkward feature of the fine structure theory of the $J_\alpha$'s is that
special parameters are required to make good sense of the notion of
``$\Sigma_n$ Skolem hull'' for $n>1.$  The source of the problem is that
parameters are needed to uniform $\Sigma_n$  relations when $n>1.$ 

The purpose of this article is to indicate how a reformulation of Jensen's
$\Sigma^*$ theory (developed for the study of core models) can be used to
provide a more satisfactory treatment of uniformization, hulls and Skolem
functions for the $J_\alpha$'s. Then we use this approach to fine structure to
formulate a principle intended to capture the combinatorial content of the
axiom $V=L.$ 

\vskip10pt

\flushpar
{\bf Section One \,  Fine Structure Revisited}

We begin with 
a simplified definition of the $J$-hierarchy. Inductively we define
$\widetilde{J}_\alpha,$ $\alpha\in\ORD$ (and then $J_\alpha=\widetilde
{J}_{\omega\alpha }):$ $\widetilde{J}_n=V_n$ for $n\le\omega.$  Suppose
$\widetilde{J}_\lambda$ is defined for a limit $\lambda$  and let
$W^\lambda_n(e,x)$ be a canonical universal $\Sigma_n(\widetilde{J}_\lambda)$
predicate (also defined inductively). For $e\in\widetilde{J}_\lambda$ let
$X^\lambda_1(e)=\{x|W^\lambda_1(e,x)\}$  and for $n\ge 1,$
$X^\lambda_{n+1}(e)=$ $\{X_n^\lambda(\bar e)|W^\lambda_{n+1}(e,\bar e)\}.$ 
Then $\widetilde{J}_{\lambda+n}=\{X^\lambda_n(e)|e\in\widetilde{J}_\lambda\}.$
For all limit
$\lambda,\widetilde{J}_\lambda=\bigcup
\{\widetilde{J}_\delta|\delta<\lambda\}.$  It is 
straightforward to verify that the $\widetilde{J}_\lambda,\lambda$ limit
behave like, and in fact equal, the usual $J_\alpha$'s.

Let $M$ denote some $J_\alpha,\alpha>0.$  (More generally, our theory applies
to ``acceptable $J$-models''.) We make the following definitions, inductively.

1) \  A $\Sigma_1^*$ formula is just a $\Sigma_1$  formula. A predicate is  
$\underline{\Sigma_1^*}$  ($\Sigma _1^*,$ respectively) if it is  definable by
a $\Sigma_1^*$ formula with (without, respectively) parameters.
$\rho^M_1=\Sigma_1^*$ projectum of $M=$ least $\rho$ s\.t\. there is a
$\underline{\Sigma_1^*}$ subset of $\omega\rho$ not in $M.$
$H^M_1=H^M_{\omega\rho_{1}^M}=$ sets $x$ in $M$ s\.t\. $M$-card (transitive
closure $(x)$) $<\omega\rho^M_1.$  For any $x\in M,$ $M_1(x)=$ First
reduct of $M$ relative to $x=\langle H^M_1, A_1(x)\rangle$  where
$A_1(x)\subseteq H_1^M$ codes the $\Sigma_1^*$ theory of $M$  with parameters
from $H_1^M\cup\{x\}$ in the natural way: $A_1(x)=\{\langle y,n\rangle|$ the
$n^{\text{th}}$ $\Sigma^*_1$ formula is true at $\langle y,x\rangle,$ $y\in
H^M_1\}.$  A good $\Sigma^*_1$ function is just a $\Sigma_1$ function and for
any $X\subseteq M$ the $\Sigma_1^*$ hull $(X)$ is just the $\Sigma_1$ hull of
$X.$ 

2) \  For $n\ge 1,$ a $\Sigma^*_{n+1}$  formula is one of the form
$\varphi(x)\longleftrightarrow M_n(x)\models\psi,$ where $\psi$ is $\Sigma_1.$
A predicate is $\underline{\Sigma^*_{n+1}}$ ($\Sigma^*_{n+1},$ respectively)
if it is defined by a $\Sigma^*_{n+1}$  formula with (without, respectively)
parameters. $\rho^M_{n+1}=\Sigma^*_{n+1}$ projectum of $M=$ least $\rho$  such
that there is a $\underline{\Sigma^*_{n+1}}$ subset of $\omega\rho$ not in
$M.$  $H^M_{n+1}=H^M_{\omega\rho^{M}_{n+1}}=$ sets $x$ in $M$  s\.t\.
$M$-card (transitive closure $(x))<\omega\rho^M_{n+1}.$ For any $x\in M,$ 
$M_{n+1}(x)=(n+1)$ s\.t\. reduct of $M$  relative to $x=\langle H^M_{n+1},
A_{n+1}(x)\rangle$  where $A_{n+1}(x)\subseteq H^M_{n+1}$ codes the
$\Sigma^*_{n+1}$ theory of $M$  with parameters from $H^M_{n+1}\cup\{x\}$ 
in the natural way: \ $A_{n+1}(x)=\{\langle y,m\rangle|$ the $m^{\text{th}}$ 
$\Sigma^*_{n+1}$ formula is true at $\langle y,x\rangle,y\in H^M_{n+1}\}.$  A
good $\Sigma^*_{n+1}$ 
function $f$ is a function whose graph is $\Sigma^*_{n+1}$ 
with the additional property that for $x\in\Dom(f),$ $f(x)\in\Sigma^*_n$ hull
$(H^M_n\cup\{x\}).$  The $\Sigma^*_{n+1}$ hull $(X)$ for $X\subseteq M$ is the
closure of $X$  under good $\Sigma^*_{n+1}$  functions.

\vskip10pt

\flushpar
{\bf Facts.} \  (a) \  $\varphi,\psi\Sigma^*_n$ formulas
$\longrightarrow\varphi\vee\psi,\varphi\wedge\psi$  are  $\Sigma^*_n$
formulas.

(b)  \  $\phi\Sigma^*_n$ or $\prod^*_n$ ($=$ negation of
$\Sigma^*_n)\longrightarrow\phi$ is $\Sigma^*_{n+1}.$

(c) \  $Y\subseteq \Sigma ^*_n$ hull $(X)\longrightarrow\Sigma^*_n$ hull
$(Y)\subseteq\Sigma^*_n$ hull $(X).$

(d)  \  $f$  good $\Sigma^*_n$ function $\longrightarrow f$ good
$\Sigma^*_{n+1}$ function.

(e) \  $\Sigma^*_n$ hull $(X)\subseteq\Sigma^*_{n+1}$ hull $(X).$

(f) \  There is a $\Sigma^*_n$ relation $W(e,x)$ s\.t\. if $S(x)$ is
$\Sigma^*_n$  then for some $e\in\omega,$ $S(x)\longleftrightarrow W(e,x)$ for
all $x.$

(g) \  The structure $M_n(x)=\langle H^M_n, A_n(x)\rangle$ is amenable.

(h) \  $H^M_n=J^{A_{n}}_{\omega\rho^{M}_{n}}$ where $A_n=A_n(0).$

(i) \  Suppose $H\subseteq M$ is closed under good $\Sigma^*_n$ functions and
$\pi:\ \overline{M}\longrightarrow M,\overline{M}$ transitive, 
Range$(\pi)=H.$ Then $\pi$ preserves $\Sigma^*_n$ formulas: for
$\Sigma^*_n\phi$ and $x\in\overline{M},$ $\oM\models\phi(x)\longleftrightarrow
M\models\phi(\pi(x)).$ 

\flushpar
{\bf Proof of (i).} \   Note that $H\cap M_{n-1}(\pi(x))$ is
$\Sigma_1$-elementary in $M_{n-1}(\pi(x)).$  And $\pi^{-1}[H\cap
M_{n-1}(\pi(x))]=\langle J^A_{\omega\rho },A(x)\rangle$   
for some $\rho,A,A(x).$  But (by induction on $n)$  $A=A^M_{n-1}\cap
J^A_{\omega\rho },$  $A(x)=A_{n-1}(x)^M\cap J^A_{\omega\rho }$  and
$\rho=\rho^M_{n-1}.$  \hfill{$\dashv$ }

\proclaim{Theorem 1} By induction on $n>0:$ 

1) \  If $\phi(x,y)$ is $\Sigma^*_n$ then $\exists y\in\Sigma^*_{n-1}$ hull
$(H^M_{n-1}\cup\{x\})\phi(x,y)$ is also $\Sigma^*_n.$

2) \  If $\phi(x_1\cdots x_k)$ is $\Sigma^*_m,m\ge n$ and $f_1(x),\cdots,
f_k(x)$ are good $\Sigma^*_n$ functions, then $\phi(f_1(x)\cdots f_k(x))$ is
$\Sigma^*_m.$

3)  \  The domain of a good $\Sigma^*_n$  function is $\Sigma^*_n.$

4) \  Good $\Sigma^*_n$ functions are closed under composition.

5) \  ($\Sigma^*_n$ Uniformization) If $R(x,y)$ is $\Sigma^*_n$ then there is a
good $\sg$ function $f(x)$ s\.t\. $x\in\Dom(f)\longleftrightarrow$ $\exists
y\in\sgn$ hull $(H^M_{n-1}\cup\{x\})R(x,y)\longleftrightarrow R(x,f(x)).$

6) \   There is a good $\sg$ function $h_n(e,x)$ s\.t\. for each $x,$ $\sg$
hull $(\{x\})=\{h_n(e,x)|e\in\omega\}.$  
\endproclaim

\demo{Proof}  The base case $n=1$ is easy (take $\so$ hull $(X)=M$ for all
$X).$   Now we prove it for $n>1,$ assuming the result for smaller $n.$

1) \  Write $\exists y\in\sgn$ hull $(H^M_{n-1}\cup\{x\})\phi(x,y)$ as 
$\exists\bar y\in H^M_{n-1}\phi(x,h_{n-1}(e,\langle x,\bar y\rangle))$ using
6) for $n-1.$  Since $h_{n-1}$ is good $\sgn$ we can apply 2) for $n-1$ to
conclude that $\phi(x,h_{n-1}(e,\langle x,\bar y\rangle))$ is $\sg.$ Since the
quantifiers $\exists e\exists\bar y\in H^M_{n-1}$ range over $H^M_{n-1}$ they
preserve $\sg$-ness.

2)  \  $\phi(f_1(X)\cdots f_k(x))\longleftrightarrow\exists x_1\cdots
x_k\in\sgn$ hull $(H^M_{n-1}\cup\{x\})$  $[x_i=f_i(x)$ for $1\le i\le
k\wedge\phi(x_1\cdots x_k)].$  If $m=n$ then this is $\sg$ by 1). If $m>n$
then reason as follows: the result for $m=n$  implies that $A_n(\langle
f_1(x)\cdots f_k(x)\rangle)$  is $\Delta_1$ over $M_{n+1}(x).$  Thus
$A_{m-1}(\langle f_1(x)\cdots f_k(x)\rangle)$ is $\Delta_1$ over $M_{m-1}(x).$
So as $\phi$ is $\Sigma^*_m$ we get that $\phi(f_1(x)\cdots f_k(x))$ is also
$\Sigma_1$  over $M_{m-1}(x),$  hence $\Sigma^*_m.$

3) \  If $f(x)$ is good $\Sigma^*_n$ then dom$(f)=\{x|\exists
y\in\Sigma^*_{n-1}$ hull of $H^M_{n-1}\cup\{x\}(y=f(x))\}$ is $\sg$ by  1).

4)  \  If $f,g$ are good $\sg$ then the graph of $f\circ g$ is $\sg$  by 2).
And $f\circ g(x)\in\sgn$ hull$(H^M_{n-1}\cup\{x\})$ since the latter hull
contains $g(x),f$ is good $\sg$ and Fact c) holds.

5)  \  Using 6) for $n-1,$ let $\overline{R}(x,\bar y)\longleftrightarrow
R(x,h_{n-1}(\bar y))\wedge\bar y\in H^M_{n-1}.$  Then $\overline{R}$ is $\sg$
by 2) for $n-1$ and using $\Sigma_1$ uniformization on $(n-1)$ s\.t\. reducts
we can define a good $\sg$ function $\bar f$ s\.t\. $\overline{R}(x,\bar
f(x))\longleftrightarrow\exists\bar y\in H^M_{n-1}\overline{R}(x,\bar y).$
Let $f(x)=h_{n-1}(\bar f(x)).$ Then $f$  is good $\sg$ by 4).

6) \  Let $W$ be universal $\sg$ as in Fact f). By 5) there is a good $\sg$
$g(e,x)$ s\.t\. $\exists y\in\sgn$ hull$(H^M_{n-1}\cup\{x\})$ $W(e,\langle
x,y\rangle)\longleftrightarrow W(e,\langle x,g(e,x)\rangle)$  (and $g(e,x)$
defined $\longrightarrow W(e,\langle x,g(e,x)\rangle)).$  Let
$h_n(e,x)=g(e,x).$  If $y\in\sg$ hull $(\{x\})$ then for some $e, W(e,\langle
x,y'\rangle)\longleftrightarrow y'=y$ so $y=h_n(e,x).$  Clearly
$h_n(e,x)\in\sg$ hull $(\{x\})$ since $h_n$ is good $\sg.$ \hfill{$\dashv$} 

\enddemo

\vskip20pt

\flushpar
{\bf Section Two \  The Combinatorial Content of $V=L$ }

In this section we provide an axiomatic treatment of the $\Sigma^*$  theory
introduced in Section One.  When establishing combinatorial principles in
$L[R], R$ a real, one makes use of a {\it standard Skolem system} for $R$
(defined below), of which the system of canonical $\sg$ Skolem functions for
the $J^R_\alpha$'s constitutes the canonical example. Our principal goal is to
provide combinatorial axioms for a system of functions which guarantee that it
is in fact a standard Skolem system for some real. These axioms can then be
used to formulate a single combinatorial principle which captures the full
power of Jensen's fine structure theory.

Some notation: \  For $\delta=\lambda+n,\lambda$ limit or $0$  and
$n\in\omega,$  $\Seq(\delta)$ denotes all finite sequences from $\lambda$
together with all finite sequences from $\delta$  of length $\le n.$  Let
$x*y$ denote the concatenation of the sequences $x,y.$  For $\lambda$ limit or
$O,\widetilde{J}^R_\lambda$ denotes $J^R_\delta$ where
$\omega\cdot\delta=\lambda.$  

A {\it standard Skolem system} for a real $R$  is a system $\vec{F}=\langle
F^\delta_n|n>0,\,\delta\in\ORD, n>1\longrightarrow\delta$ limit$\rangle$ where
$F^\delta_n$ is a partial function from $\omega\times\Seq(\delta)$ to
$\delta,$  obeying  (A) -- (E) below. For any limit
$\lambda,x\in\Seq(\lambda), n\ge 1$ let
$H^\lambda_n(x)=\{F^\lambda_n(k,x)|k\in\omega\}$  and if $\bar\lambda=$
ordertype $(H^\lambda_n(x))$ let $\pi^n_{\bar\lambda\lambda }(x):\
\bar\lambda\longrightarrow\lambda$  be the increasing enumeration of
$H^\lambda_n(x).$  We say $y\in H^\lambda_n(x),$ for $y\in\Seq(\lambda),$ if
$y^*\in H^\lambda_n(x)$ where $y^*$ is a canonical ordinal code for $y.$

{\bf (A) \  (Monotonicity)}  \  $\delta_1\le\delta_2\longrightarrow
F^{\delta_{1}}_1\subseteq F^{\delta_{2}}_1, x\in H^\lambda_1(x)\subseteq
H^\lambda_2(x)\subseteq\cdots\subseteq\lambda$ for limit
$\lambda,x\in\Seq(\lambda).$

{\bf (B) \ (Condensation)} \  Let $\pi=\pi^n_{\bar\lambda\lambda }(x).$  
Then for
$m\le n,$ and $\bar x\in\Seq(\bar\lambda),$ $\pi(F^{\bar\lambda }_m(k,\bar
x))\simeq F^\lambda_m(k,\pi(\bar x)).$  And
$\tilde\pi(F^{\bar\lambda+m}_1(k,\bar x))\simeq
F^{\lambda+m}_1(k,\tilde\pi(\bar x))$  for $\bar x\in\Seq(\bar\lambda+m),$
where $\tilde\pi$ is the extension of $\pi$  to  $\bar\lambda+m$  obtained by
sending $\bar\lambda+i$ to $\lambda+i.$

{\bf (C)  \  (Continuity)} \  For limit
$\lambda,F^\lambda_1=\bigcup\{F^\delta_1|\delta<\lambda\}.$  There is a
$p\in\Seq(\lambda)$  such that for all $x\in\Seq(\lambda)$  and $y<\lambda,$
$F^\lambda_{n+1}(x)\simeq y$ iff for some $z\in\Seq(\lambda), F^{\bar\lambda
}_{n+1}(\bar x)\simeq\bar y$  where $\bar\lambda=$ ordertype
$(H^\lambda_n(z)),$ $\pi^n_{\bar\lambda\lambda }(z)$ sends $\bar x,\bar y$ to
$x,y$  and $p\in H^\lambda_n(z).$

{\bf (D) } \    $\langle F^\delta_n|\delta<\lambda,n<\omega\rangle$  is
uniformly $\Delta_1(\widetilde{J}^R_\lambda)$  for limit $\lambda,$ in the
parameter $R.$

{\bf (E)} \  For limit $\lambda,H^\lambda_1(x)=\lambda\cap\Sigma_1$ Skolem
hull of $x$  in $\widetilde{J}_\lambda^R$ for $x\in\Seq(\lambda)$  and for
some fixed $p\in\Seq(\lambda), \bigcup\limits_n
H^\lambda_n(x)=\lambda\cap$ Skolem hull of
$x$  in $\widetilde{J}^R_\lambda$  whenever $p$  belongs to $H^\lambda_n(x)$
some $n,x\in\Seq(\lambda).$

Intuitively, $F^\lambda_n$ is a $\sg$ Skolem function for
$\widetilde{J}^R_\lambda$  and $F^{\lambda+n}_1$ is the $n^{\text{th}}$
approximation to $F^{\lambda+\omega }_1.$

\proclaim{Proposition 2} For every real $R$  there exists a standard Skolem
system for $R.$ 
\endproclaim

\demo{Proof} Let $\psi\longmapsto\psi^*_n$ be a recursive translation on
formulas so that for limit $\lambda,$
$\widetilde{J}^R_{\lambda+n}\models \psi \longleftrightarrow
\widetilde{J}^R_\lambda\models\psi^*_n$
(where $\widetilde{J}^R_\alpha$ is defined just like $\widetilde{J}_\alpha,$
but relativized to $R).$  Fix a recursive enumeration
$\langle\phi_k(v)|k\in\omega\rangle$  of $\Delta_0$ formulas with a predicate
$\underline{R}$ denoting $R$  and sole free variable $v.$  Let $<_R$ denote
the ordering of $L[R]$  given by: \  $x<_Ry$ iff
$\exists\lambda\in\Lim\cup\{0\}\exists n\in\omega$
$[y\in\widetilde{J}^R_{\lambda+n+1}-\widetilde{J}^R_{\lambda+n},$
$(x\in\widetilde{J}^R_{\lambda+n})$ or $(\lambda$ limit,
$x\in\widetilde{J}^R_{\lambda+n+1},$ $e<_Rf$ where $e,f$ are $<_R-$ least
s\.t\. $X^{\lambda,R}_{n+1}(e)=x,$ $X^{\lambda,R}_{n+1}(f)=y)$ or $(\lambda=0$
and $x<_Ly)].$   

Now define $\vec F=\langle F^\delta_n|\delta\in\ORD, n>0,
n>1\longrightarrow\delta\text{ limit}\rangle$ as follows:

(a)  \  $F^n_1(k,x)\simeq y$ iff $L^R_n\models\exists w$ s\.t\. $\langle
y,w\rangle$ is $<_R-$ least s\.t\. $\phi_k(\langle x,y,w\rangle).$

(b)  \  For $\lambda$ limit, $F^\lambda_1=\cup\{F^\delta_1|\delta<\lambda\}.$

(c)  \  For $\lambda$ limit, $n>0,$ $F^{\lambda+n}_1(k,x)\simeq y$ iff for
some $m\le n,$ $\widetilde{J}^R_{\lambda+m}\models(\exists w$ s\.t\. $\langle
y,w\rangle$  is $<_R-$ least s\.t\. $\phi_k(\langle x,y,w\rangle))$ and if
$\psi$ denotes the formula in parentheses then $\psi^*_m$ is $\sg.$

(d)  \  For $\lambda$ limit, $n>1,$ $F^\lambda_n$ is the canonical $\sg$
Skolem function for $\widetilde{J}^R_\lambda$ (restricted to
$\omega\times\Seq(\lambda))$ as in 6) of Theorem 1.

The verification that $\vec F$ is a standard Skolem system for $R$  is
straightforward as Condensation is guaranteed by (c) above and (C), (E) are
satisfied by letting $p$  be the full standard parameter for
$\widetilde{J}^R_\lambda.$  \hfill{$\dashv$} 
\enddemo

An {\it abstract Skolem system} is a system $\vec F$ obeying properties (A),
(B), (C) from the definition of standard Skolem system. We would like to prove
that every abstract Skolem system is a standard Skolem system for some real.
However standard systems share one further property which we must also impose:

\flushpar
(Stability) \  For $\lambda$ limit, $x\in\Seq(\lambda)$ let $\pi:\
\bar\lambda\longrightarrow\lambda$  be the increasing enumeration of
$H^\lambda_1(x).$  Then $\pi$ extends uniquely to a $\Sigma_1$-elementary
embedding of $\langle\widetilde{J}_{\bar\lambda }^{\vec F},\vec
F\restriction\bar\lambda\rangle$  into $\langle\widetilde{J}^{\vec
F}_\lambda,\vec F\restriction\lambda\rangle.$ Also for $\lambda$ limit there
is $p\in\Seq(\lambda)$ such that for all $x\in\Seq(\lambda),$ if $\pi:\
\bar\lambda\longrightarrow\lambda$ is the increasing enumeration of
$H^\lambda(x)=\bigcup\limits_n H^\lambda_n(x)$  and $p\in H^\lambda(x)$ then
$\pi$ extends uniquely to an elementary embedding of
$\langle\widetilde{J}^{\vec F}_{\bar\lambda },\vec
F\restriction\bar\lambda\rangle$  into $\langle\widetilde{J}^{\vec F}_\lambda,
\vec F\restriction\lambda\rangle.$   

Though stability is not combinatorial we
shall see that any abstract Skolem system can be made stable without changing
its ``cofinality function''. This fact will enable us to formulate
combinatorial principles which are universal for principles which depend only
on cofinality.

\proclaim{Theorem 3} The following are equivalent:

(a)  \  $\vec F$ is a stable, abstract Skolem system.

(b) \  $\vec F$ is a standard Skolem system in a $CCC$  forcing extension of
$V.$ 

\endproclaim

\flushpar
Note that (b) $\longrightarrow$  (a) follows easily, using the absoluteness of
the concept of 
stability. We now develop the forcing required to prove (a) $\longrightarrow$
(b). 

Fix a stable, abstract Skolem system $\vec F$  and let $M$  denote $L[\vec
F],$ $M_\lambda=\langle\widetilde{J}^{\vec F}_\lambda,\vec
F\restriction\lambda\rangle$  for limit $\lambda.$   The desired forcing
${\Cal{P}}$ is a $CCC$  forcing of size $\omega_1$ in $M.$  It is designed so
as to produce a generic real $R$  which codes $\vec F\restriction\omega_1$ via
a careful almost disjoint coding. We will demonstrate that $R$  in fact codes
all of $\vec F$ using condensation properties of $\vec F.$ 

We begin our description of ${\Cal{P}}.$  A limit ordinal $\lambda$ is {\it
small} if for some $x\in\Seq(\lambda)$ and some $n,H^\lambda_n(x)=\lambda.$
Let $n(\lambda)$  be the least $n$  s\.t\. such an $x$  exists and let
$p^\lambda$  be the least $p\in\Seq(\lambda)$ s\.t\.
$H^\lambda_{n(\lambda)}(p)=\lambda.$  We now define a canonical bijection
$\bar f_\lambda:\ \lambda\longrightarrow\omega.$  First let
$g:\lambda \longrightarrow\omega$ be defined by $g(\delta)=$ least $k$  s\.t\.
$\delta=F^\lambda_{n(\lambda)}(k,p^\lambda).$  Then $\bar f_\lambda(\delta)=m$
if $g(\delta)$ is the $m^{\text{th}}$ element of Range$(g)$ under $<$ on
$\omega.$ Now let $f_\lambda:\omega\longrightarrow M_\lambda$  be
$g^*\circ\bar f^{-1}_\lambda$  where $g^*:\lambda\longrightarrow M_\lambda$
is a canonical $\underset{\sim}\to\Delta_1(M_\lambda)$ bijection. Now choose
$A_\lambda\subseteq\omega$ to code $M_\lambda$ using $f_\lambda$ and let
$b_{\lambda+n(\lambda)}$  be a function from $\omega$ to $\omega$ which is
$\Delta_{n(\lambda)+1}\langle M_\omega,A_\lambda\rangle$ yet eventually
dominates each function from $\omega$ to $\omega$ which is
$\Delta_{n(\lambda)}\langle M_\omega ,A_\lambda \rangle.$ Also require
that Range$(b_{\lambda+n(\lambda)})\subseteq_*$ Range$(b_{\bar\lambda+n})$ for
all 
$\bar\lambda<\lambda,n<\omega$  where we have (inductively) defined
$b_{\bar\lambda+n}.$ ($\subseteq_*$ denotes inclusion except for a finite
set.) 

We also define $b_{\lambda+n}$  for $n=n(\lambda)+m,m>0.$  For this purpose
define $\overline{F}_1^{\lambda+n}(k,\bar x)\simeq\bar y$ to mean
$F_1^{\lambda+n}(k,x)\simeq y$  where $x(i)=\lambda+\bar x(i)$ if $\bar
x(i)<n,$ $\bar x(i)=n+x(i)$  otherwise (similarly for $y).$  Let
$A_{\lambda+m}\subseteq\omega$ code.
$$\langle
M_\lambda,\overline{F}_1^{\lambda+n(\lambda)},\cdots,\overline{F}_1^{ 
\lambda+n(\lambda)+m-1}, F^\lambda_{n(\lambda)},\cdots,
F^\lambda_{n(\lambda)+m -1}\rangle $$
using $f_\lambda$  and let $b_{\lambda+n(\lambda)+m}$ be a function from
$\omega$  to $\omega$ which is $\Delta_{n(\lambda)+m+1}\langle
M_\omega,A_{\lambda+m}\rangle$ 
yet eventually dominates
$\Delta_{n(\lambda)+m}\langle M_\omega,A_{\lambda+m}\rangle$    functions.
Also require that Range$(b_{\lambda+n(\lambda)+m})\subseteq_*$
Range$(b_{\lambda+n(\lambda)+m-1}).$  We use the $b_{\lambda+n},n\ge
n(\lambda)$  to facilitate the desired almost disjoint coding.

An {\it index} is a tuple of one of the forms $\langle\lambda+n,1,k,\bar
x,\bar y\rangle,$ $\langle\lambda,n,k,\bar x,\bar y\rangle$ where $\lambda$ is
small, $n\ge n(\lambda)$ and $\overline{F}^{\lambda+n}_1(k,\bar x)\simeq\bar
y,$ $F^\lambda_n(k,\bar x)\simeq\bar y,$  respectively. Let $\langle
Z_e|e\in\omega\rangle$ be a recursive partition of $\omega-\{0\}$ into
infinite pieces. For each index $x$  we define a ``code'' $b_x$  as follows: \
 If $x=\langle\lambda+n,1,k,\bar x,\bar y\rangle,$ 
$\langle\lambda,n,k,\bar x,\bar y\rangle$ then $b_x=b_{\lambda+n}\restriction
Z_e$ where $f_\lambda(e)=\langle n,1,k,\bar x,\bar y\rangle,$ 
$\langle 0,n,k,\bar x,\bar y\rangle,$  respectively. A {\it restraint} is a
function of the form $b_x,x$ an index. We sometimes view $b_x$ as a subset of
$\omega$ by identifying it with $\{\langle n,m\rangle|b_x(n)=m\},$
$\langle\cdot \, , \,\cdot \rangle$ a recursive pairing on $\omega.$ 

A condition in ${\Cal{P}}$ is $p=\langle s,\bar s\rangle$  where
$s:|s|\longrightarrow 2,$ $|s|\in\omega,$ $\bar s$ is a finite set of
restraints and when $i=\langle m,k,x,y\rangle<|s|$  then
$s(i)=1\longleftrightarrow F^m_1(k,x)\simeq y.$  Extension is defined by: \
$(s,\bar s)\le(t,\bar t)$ iff $s\supseteq t,$ $\bar s\supseteq\bar t$ and
$s(i)=1\longrightarrow t 
(i)=1$  or $i\notin\bigcup\bar t.$ (Recall that we can
think of $b_x\in\bar t$ as a subset of $\omega.)$ 

This is a $CCC$  forcing and a generic $G$  is uniquely determined by the real
$R=\bigcup\{s|(s,\bar s)\in G$ for some $\bar s\}.$  Fix such a real $R.$  

\proclaim{Lemma 4} $\langle F^\delta_n|\delta<\lambda,n<\omega\rangle$ is
uniformly $\Delta_1(\widetilde{J}^R_\lambda)$ for limit $\lambda,$ in the
parameter $R.$ 
\endproclaim

\demo{Proof} By induction we define $F^\lambda_n,$ $F^{\lambda+n}_1$ for
$\lambda$ limit or $0,$ $n\in\omega.$ If $\lambda=0$  then $F^n_1$ can be
defined directly from $R$  by the restriction we placed on $s$  for conditions
$(s,\bar s).$  For $\lambda$ limit, $F^\lambda_1$ is defined by induction and
Continuity. Also, induction and Continuity enable us to define $F^\lambda_n,$
$F^{\lambda+n}_1$  provided $n\le n(\lambda)\ne 1$ or $n(\lambda)$ is not
defined. Thus if $\lambda$ is not small we're done and otherwise we can define
$f_\lambda,b_{\lambda+n},$ by induction. Let $f_\lambda(e)=\langle n,1,k,\bar
x,\bar y\rangle.$  Then $\overline{F}_1^{\lambda+n}(k,\bar x)\simeq\bar y$ iff
$\langle\lambda+n,1,k,\bar x,\bar y\rangle$ is an index iff $R$ is almost
disjoint from $b_{\lambda+n}\restriction Z_e.$  The definition of
$F^\lambda_n$ is similar, using $\langle 0,n,k,\bar x,\bar y\rangle.$ 
\hfill{$\dashv$} 
\enddemo

Our next goal is to establish a strong statement of the definability of the
forcing relation for ${\Cal{P}}.$  For any infinite ordinal $\delta$  we let
${\Cal{P}}(\delta)$ denote those conditions in ${\Cal{P}}$ involving
restraints with indices $\langle\lambda+n,1,k,\bar x,\bar y\rangle,$
$\langle\lambda,n,k,\bar x,\bar y\rangle$ where $\lambda+n<\delta.$  For $p\in
{\Cal{P}}$ we let $p\restriction\delta$ be obtained from $p$  by discarding
all restraints which are not of the above form.

\proclaim{Lemma 5} \  (Persistence) let $\lambda$ be small and for $p\in
{\Cal{P}}(\lambda+\omega)$  let $p^*$ be obtained by replacing each of its
restraints of the form $b_x,x=\langle\lambda+n,1,k,\bar x,\bar y\rangle,$
$\langle\lambda,n,k,\bar x,\bar y\rangle$  by $\langle n,1,k,\bar x,\bar
y\rangle,$ $\langle n,k,\bar x,\bar y\rangle,$  respectively. (Then $p^*\in
M_\lambda.)$  Suppose $W\subseteq {\Cal{P}}(\lambda+n(\lambda)+m)$ and
$W^*=f_\lambda^{-1}[\{p^*|p\in W\}]$ is $\Sigma_{n(\lambda)+m}$  over $\langle
M_\omega,A_{\lambda+m}\rangle.$  Then
$D=\{p\in{\Cal{P}}(\lambda+n(\lambda)+m)|\exists q\in W(p\le q)$ or $\
\forall\ q\le p(q\notin W)\}$ is predense on ${\Cal{P}}.$  
\endproclaim

\demo{Proof} Given $p\in {\Cal{P}}$ we must find $q\le p$ such that
$q\restriction\lambda+n(\lambda)+m$ belongs to $D.$  Write $p=(s,\bar
s\cup\bar t)$ where $p\restriction\lambda+n(\lambda)+m=(s\bar s),$ $\bar s\cap
\bar t=\emptyset.$ For each $n$  let $s_n$  extend $s$ by assigning $\langle
m_0,m_1\rangle$ to $0$  whenever $\langle m_0,m_1\rangle\notin\Dom(s)$  and
$m_0\le m_1\le n.$  (We intend that $n\longmapsto s_n$ is recursive.) If
$(s_n,\bar s)$ belongs to $D$  for some $n$  then we are done since $(s_n,\bar
s\cup\bar t)$  extends $p.$  If not then we can define a
$\Sigma_{n(\lambda)+m}$  over $\langle M_\omega,A_{\lambda+m}\rangle$
function $n\longmapsto t_n$  so that for some $\bar t_n,$ $(t_n,\bar
t_n)\le(s_n,\bar s),$ $(t_n,\bar t_n)\in W,$ using the fact that
$A_{\lambda+m}$ codes $\langle
M_\lambda,\overline{F}^{\lambda+n(\lambda)}_1,\cdot,\overline{F}^{\lambda+n(\lambda)+m-1}_1\rangle$  and
hence ``codes'' ${\Cal{P}}(\lambda+n(\lambda)+m).$  Then $f(m+1)=$ length
$(t_{f(m)}),f(0)=0$ defines a $\Sigma_{n(\lambda)+m}$  over $\langle
M_\omega,A_{\lambda+m}\rangle$  function and every such function is eventually
dominated by the function $b_{\lambda+n(\lambda)+m}.$  Thus there must be
infinitely many $\ell$  such that $[f(\ell), f(\ell+1)]$ is disjoint from
Range$(b_{\lambda+n(\lambda)+m}).$  As Range$(b)\subseteq_*$
Range$(b_{\lambda+n(\lambda)+m})$ for all $b\in\bar t$  it follows that for
some $\ell, [f(\ell), f(\ell+1)]$ is disjoint from $\cup\{\text{Range}(b)|b\in
\bar t\}.$ But then $(t_{f(\ell)},\bar t_{f(\ell)}\cup\bar t)=q\le q$  and
$q\restriction\lambda+n(\lambda)+m$ belongs to $W\subseteq D.$
\hfill{$\dashv$ } 
\enddemo

\proclaim{Corollary 6} The forcing relation
$\{(p,\phi)|p\in{\Cal{P}}(\lambda)$  and $p\Vdash\phi$ in ${\Cal{P}}(\lambda)$
where $\phi$ is a ranked sentence in $M_\lambda\}$ is $\Sigma_1$ over
$M_\lambda,$ for limit $\lambda.$ 
\endproclaim

\demo{Proof} By induction on $\lambda.$  Note that if
$\bar\lambda<\lambda,\bar\lambda$ limit then for
$p\in{\Cal{P}}(\bar\lambda),\phi$ ranked in $M_{\bar\lambda }$ we have
$p\Vdash\phi$ in ${\Cal{P}}(\bar\lambda)$ iff $p\Vdash\phi$ in
${\Cal{P}}(\lambda).$  The reason is that by Lemma 5, every
${\Cal{P}}(\lambda)$-generic is ${\Cal{P}}(\bar\lambda)$-generic for ranked
sentences, since by induction the ${\Cal{P}}(\bar\lambda)$ forcing relation
for ranked sentences is $\Sigma_1$ over $M_\lambda.$ 

Thus we are done by induction if $\lambda$ is a limit of limit ordinals. Now
suppose that we wish to establish the Corollary for $\lambda+\omega.$  We may
assume that $\lambda$ is small as otherwise ${\Cal{P}}(\lambda+\omega)$ is a
set forcing in $M_{\lambda+\omega }.$  Now any ranked sentence $\phi$  in
$M_\lambda$ is equivalent to a $\Sigma_{n(\lambda)+M}$ statement about
$M_\lambda [\underline{R}]$ for some $m(\underline{R}$ denoting the generic
real). But then by Lemma 5, $p\Vdash\phi$ in ${\Cal{P}}(\lambda+\omega)$ iff
$p\Vdash\phi$ in ${\Cal{P}}(\lambda+n(\lambda)+m)$ for $p\in
{\Cal{P}}(\lambda+n(\lambda)+m).$  As the latter is $\Sigma_1$-definable over
$M_{\lambda+\omega },$ we are done. \hfill{$\dashv$ }
\enddemo

\proclaim{Corollary 7} Suppose $\lambda$ is small and
$W\subseteq{\Cal{P}}(\lambda)$ is $\Sigma_{n(\lambda)}$ over $M_\lambda.$ Let
$D=\{p\in {\Cal{P}}(\lambda)|\exists q\in W(p\le q)$ or $\ \forall\ q\le
p(q\notin W)\}.$  Then $D$ is predense on ${\Cal{P}}.$  
\endproclaim

\demo{Proof} Let $m=0$ in Lemma 5. \hfill{$\dashv$ }

\enddemo

Now we are prepared to finish the proof of the Characterization Theorem. Note
that the only remaining condition to verify in showing that $\vec F$ is a
Standard Skolem system is condition (E), where stability is used.

\proclaim{Lemma 8} For $\lambda$ limit, $x\in\Seq(\lambda),$
$H^\lambda_1(x)=\lambda\cap\Sigma_1$ Skolem hull of $x$  in
$\widetilde{J}^R_\lambda.$  For $\lambda$ limit there is $p\in\Seq(\lambda)$
s\.t\. for all $x\in\Seq(\lambda),$ $H^\lambda(x)=\bigcup\limits_n
H^\lambda_n(x)=\lambda\cap$ Skolem hull of $x$ in $\widetilde{J}^R_\lambda$
whenever $p\in H^\lambda(x).$ 
\endproclaim

\demo{Proof}  We begin with the first statement. The inclusion
$H^\lambda_1(x)\subseteq\Sigma_1$ Skolem hull of $x$  in
$\widetilde{J}^R_\lambda$ follows from Lemma 4 and Continuity. To prove the
converse we make a definition: $R$  is $\underline{\Sigma_n-\text{generic}}$
for ${\Cal{P}}(\lambda)$  if for any $\Sigma_n(M_\lambda)$
$W\subseteq{\Cal{P}}(\lambda)$ there exists $p\in G\cap {\Cal{P}}(\lambda),$
$G$  denoting the generic determined by $R,$  such that either $p$  extends a
condition in $W$  or $p$  has no extension in $W.$  By Corollary 7, if
$\lambda$  is small then $R$  is $\Sigma_{n(\lambda)}$-generic for
${\Cal{P}}(\lambda).$  

Suppose  $\phi(x,y)$ is a $\Sigma_1$ formula with parameter $x.$  Let $\pi:\
\bar\lambda\longrightarrow\lambda$  be the increasing enumeration of
$H^\lambda_1(x)$  and let $\pi(\bar x)=x.$  By Corollary 6 the forcing
relation for ${\Cal{P}}(\bar\lambda)$ is $\Sigma_1(M_{\bar\lambda })$ is
$\Sigma_1(M_{\bar\lambda })$ for ranked sentences. Since $R$ is
$\Sigma_1$-generic for ${\Cal{P}}(\bar\lambda)$ there is $p\in
G\cap{\Cal{P}}(\bar\lambda)$  s\.t\.either $p\Vdash\phi(\bar x,\bar y)$ in
${\Cal{P}}(\bar\lambda)$ for some $\bar y$  or $p\Vdash\neg\exists\bar
y\phi(\bar x,\bar y)$ in ${\Cal{P}}(\bar\lambda).$  Since $\vec F$ is stable
we have that $p\Vdash\neg\exists y\phi(x,y)$ in ${\Cal{P}}(\lambda)$  or
$p\Vdash\phi(x,y)$  where $y=\pi(\bar y).$  (Note that $\pi$  extends to a
$\Sigma_1$-elementary embedding $\tilde\pi:\ M_{\bar\lambda }\longrightarrow
M_\lambda$  such that $\tilde\pi(p)=p.)$  If $\lambda$ is small then $R$  is
$\Sigma_1$-generic for ${\Cal{P}}(\lambda)$ and thus we have shown that
$\lambda\cap\Sigma_1$ Skolem hull of $x$  in $\widetilde{J}^R_\lambda$ is
contained in $H^\lambda_1(x).$  But the above shows that if $R$ is
$\Sigma_1$-generic for ${\Cal{P}}(\lambda)$ for all small $\lambda$ then $R$
is $\Sigma_1$-generic for all $\lambda.$  So we're done.

To prove the second statement, choose $p$  to witness stability for $\vec F.$
The direction $H^\lambda(x)\subseteq$ Skolem hull of $x$  in
$\widetilde{J}^R_\lambda$ follows again from Lemma 4. For the converse, handle
each formula $\psi(x,y)$ as in the $\Sigma_1$ case, using stability and the
assumption that $p\in H^\lambda(x).$  \hfill{$\dashv$ }
\enddemo

This completes the proof of Theorem 3.

\vskip10pt

\flushpar
{\bf Universal Combinatorial Principles.}

Inherent in any abstract Skolem system $\vec F$ is its cofinality function
$\cof^{\vec F}$ defined at limit ordinals $\lambda$ as follows: $\cof^{\vec
F}(\lambda)=$ least ordertype of an unbounded subset of $\lambda$ of the form
$H^\delta_n(\gamma\cup\{p\})=\bigcup\{H^\delta_n(x*p)|x\in\Seq(\gamma)\}$ for
some $\delta\ge\lambda,n\ge 1,$ $\gamma\le\lambda,$ $p\in\Seq(\delta).$  For
any inner model $M$  let $\cof^M$ be the cofinality function of $M.$  And
$\cof=\cof^V.$ 

\proclaim{Lemma 9} Suppose $\vec F$  is an abstract Skolem system. Then there
exists a stable abstract Skolem system $\vec G$ such that $\cof^{\vec
G}=\cof^{L[\vec F]}.$ 
\endproclaim

\demo{Proof} Let $\vec G$ be obtained from $\vec F$ just as in the proof of
Proposition 2, with $R$  replaced by $\vec F.$  Then $\vec G$  is stable.
Since $\vec G$ codes $L[\vec F],$ $\cof^{\vec G}(\lambda)\le\cof^{L[\vec
F]}(\lambda)$  all $\lambda.$  But $\vec G$ is $\langle L[\vec F],\vec
F\rangle$-definable, so $\cof^{\vec G}=\cof^{L[\vec F]}.$ \hfill{$\dashv$ }
\enddemo

We now state our Universal Combinatorial Principle $P.$ 

\proclaim{Principle P} There is an Abstract Skolem System $\vec F$ such that
$\cof^{\vec F}=\cof.$ 

We show that $P$ implies all ``fine-structural principles'' for $L.$ 
\endproclaim

\proclaim{Definition} A {\bf fine-structural principle} is a statement of the
form $\exists {\Cal{A}}\psi({\Cal{A}}),$ where ${\Cal{A}}$ denotes a class and
$\psi$  is first-order, such that:

(a)  \  For every real $R$  and every Standard Skolem System $\vec F$ for
$R,L[R]\models\psi({\Cal{A}})$ for some ${\Cal{A}}$ which is definable over
$\langle L[\vec F],\epsilon,\vec F\rangle.$

(b)  \  If $M,N$  are inner models of $ZFC,$ ${\Cal{A}}$ is amenable to both
$M,N,\cof^M=\cof^N$ and $\langle M,{\Cal{A}}\rangle\models\psi({\Cal{A}})$
then $\langle N,{\Cal{A}}\rangle\models\psi({\Cal{A}}).$ 
\endproclaim

\proclaim{Theorem 10} $P$ implies all fine-structural principles.
\endproclaim

\demo{Proof} Suppose $M\models P$ with witness $\vec F$ and let $\phi$  be
fine-structural. Then $\cof^{\vec F}=\cof^M=\cof^{L[\vec F]},$ since $L[\vec
F]\subseteq M.$  By Lemma 9 there is $\vec G$ amenable to $M$  such that
$\cof^{\vec G}=\cof^M$  and $\vec G$ is stable. By the Characterization
Theorem there is a (generic) real $R$  such that $\vec G$ is a Standard Skolem
System for $R$  and hence $L[R]\models\phi$ with witness ${\Cal{A}}$ definable
over $\langle L[\vec G],\epsilon,\vec G\rangle.$  Then ${\Cal{A}}$ is amenable
to $M$  and $\cof^M=\cof^{\vec G}=\cof^{L[R]},$ so $M\models\phi.$ 
\hfill{$\dashv$ }
\enddemo

$\square$  and Morass are fine-structural but $\diamond$ is not. To obtain a
universal principle which also implies $\diamond$ we introduce a strengthening
of $P.$  

{\bf Principle $\bold{P^*}.$ } \  $V=L[\vec F]$  where $\vec F$ is an Abstract
Skolem System.

Note that $P^*\longrightarrow P,$ in view of Lemma 9. We define an {\bf L-like
principle} to be a statement $\phi$  which is true in $L[\vec F]$ whenever
$\vec F$ is a Standard Skolem System. By Lemma 9 and the Characterization
Theorem, $P^*$ implies all $L$-like principles. But unfortunately $P^*$ is not
much weaker than $V=L:$ 

\proclaim{Theorem 3.3} $P^*$ holds iff $V=L[A],A\subseteq\omega_1$ where $A$
is $L$-reshaped $(\alpha<\omega_1\longrightarrow\alpha<\omega_1$ in
$L[A\cap\alpha ]).$  
\endproclaim

\demo{Proof} Suppose $V=L[\vec F]$ for some Abstract Skolem System $\vec F.$
By Lemma 3.1 and the Characterization Theorem, we may assume that  $\vec F$ is
a Standard Skolem System for some real $R.$  Now suppose that $\alpha$  is
countable in $L[R].$  If $\alpha<\lambda$ limit,
$\widetilde{J}^R_\lambda\models\alpha$ uncountable then $\vec
F\restriction\lambda$  can be recovered inductively from $\vec
F\restriction\alpha,$  using continuity and condensation for Abstract Skolem
Systems. We can also recover $F^\lambda_n$  for all $n>0$ for such $\lambda.$
Thus if $\lambda$ is least so that $\alpha$ is countable in
$\widetilde{J}^R_{\lambda+\omega },$  we see that $\alpha$ is countable in
$L[\vec F\restriction\alpha ].$ So $\vec F\restriction\omega_1$ is
$L$-reshaped. The same argument shows that $\vec F$ is definable over $L[\vec
F\restriction\omega_1]$  so we have the desired conclusion.

For the converse note that for $L$-reshaped $A\subseteq\omega_1$  we can
define the Canonical Skolem System $\vec F^A$ for $A$  as we defined $\vec
F^R$ for reals $R,$  provided we replace the hierarchy
$\widetilde{J}^R_\delta,\delta\in\ORD$ by
$\widetilde{J}^A_\delta,\delta\in\ORD$  and we assume that for
$\lambda<\omega_1,A\cap[\lambda,\hat\lambda ]=\emptyset$ where $\hat\lambda$
is the least limit so that $\widetilde{J}^{A\restriction\lambda }_{\hat\lambda
}\models\lambda$ is countable. Then $L[\vec F^A]=L[A]$  and $\vec F^A$
satisfies the axioms for an Abstract Skolem System. (In fact $\vec F^A=\vec
F^R$  for some generic real $R$ coding $A.)$  \hfill{$\dashv$ }
\enddemo

Though $P^*$ does not therefore have models which are very far from $L,$  we
hope that its analogue in the context of core models will lead to an
interesting class of  ``$K$-like'' models.

\vskip10pt

\centerline{\bf References}

\vskip5pt

Jensen, R.B. [72] \  The Fine Structure of the Constructible Hierarchy, Annals
of Mathematical Logic.

\vskip5pt

Jensen, R.B. [89] \  Handwritten notes on the $\Sigma^*$ Theory.

\vskip5pt

Jensen, R.B. and Solovay, R.M. [68] \  Some Applications of Almost Disjoint
Sets, in {\it Mathematical Logic and the Foundations of Set Theory,}
Bar-Hillel, Editor.

\enddocument